\magnification=1200
\font\smallrm=cmr8
\def\ell{l}
\def\st{\,:\,}
\let\ppar\par
\def\reffc{1}

\nopagenumbers
\phantom{a}
\bigskip
\centerline{\bf DISCHARGING CARTWHEELS}
\bigskip\bigskip
\baselineskip=11pt

\centerline{Neil Robertson$^{*1}$\vfootnote{$^*$}{\smallrm
Research partially 
performed under a consulting agreement with Bellcore, and partially
supported by DIMACS Center, 
Rutgers University, New Brunswick, New Jersey  08903, USA.}
\vfootnote{$^1$}{\smallrm Partially supported
by NSF under Grant No. DMS-8903132 and by ONR under Grant No. 
N00014-92-J-1965.
}}
\centerline{Department of Mathematics}
\centerline{Ohio State University}
\centerline{231 W. 18th Ave.}
\centerline{Columbus, Ohio  43210, USA}
\bigskip

\centerline{Daniel P. Sanders$^{2}$\vfootnote{$^2$}{\smallrm
Partially supported by DIMACS and by ONR under Grant No. 
N00014-93-1-0325.
}}
\centerline{School of Mathematics}
\centerline{Georgia Institute of Technology}
\centerline{Atlanta, Georgia  30332, USA}
\bigskip

\centerline{P. D. Seymour}
\centerline{Bellcore}
\centerline{445 South St.}
\centerline{Morristown, New Jersey  07960, USA}
\bigskip

\centerline{and}
\bigskip

\centerline{Robin Thomas$^{*3}$\vfootnote{$^3$}{\smallrm
Partially supported
by NSF under Grant No. DMS-9303761 and by ONR under Grant No. 
N00014-93-1-0325.
}}
\centerline{School of Mathematics}
\centerline{Georgia Institute of Technology}
\centerline{Atlanta, Georgia  30332, USA}
\bigskip\bigskip\bigskip

\centerline{\bf ABSTRACT}
\bigskip

\noindent
In [{\it J.~Combin.\ Theory Ser.\ B} 70 (1997), 2-44] we gave a simplified proof
of the Four Color Theorem. The proof is computer-assisted in the sense
that for two lemmas in the article we did not give proofs, and instead
asserted that we have verified those statements using a computer. Here
we give additional details for one of those lemmas, and we include the
original computer programs and data as ``ancillary files" accompanying this
submission. 


\vfill
\noindent 30 April 1995.
 Revised  4 February 1997
\vfil
\eject

\footline{\hss\tenrm\folio\hss}
\baselineskip 22pt

\outer\def\beginsection#1\par{\vskip0pt plus.3\vsize
   \vskip0pt plus-.3\vsize\bigskip\vskip\parskip
   \message{#1}\leftline{\bf#1}\nobreak\smallskip}
\def\sqr#1#2{{\vcenter{\vbox{\hrule height.#2pt
\hbox{\vrule width.#2pt height #1pt \kern#1pt
\vrule width.#2pt}
\hrule height.#2pt}}}}
\def\square{\mathchoice\sqr56\sqr56\sqr{2.1}3\sqr{1.5}3}
\def\pf{\medskip\noindent{\sl Proof.~~}}
\def\uom{\Omega}\def\ep{\epsilon}\def\lam{\lambda}
\def\calm{{\cal M}}\def\calp{{\cal P}}\def\calr{{\cal R}}\def\calt{{\cal 
T}}

\def\claim#1#2\par{\smallbreak\noindent\rlap{\rm(#1)}\ignorespaces
\hangindent=30pt\hskip30pt{\sl#2}\smallskip}
\def\qed{\hfill$\square$\medskip}
\newcount\subsecno
\newcount\thmno
\global\subsecno=1\thmno=0
\def\thm#1. #2\par{\medbreak\global\advance\thmno by 1
\expandafter\xdef\csname#1\endcsname{(\the\subsecno.\the\thmno)}
\noindent{\bf (\the\subsecno.\the\thmno)\ }{\sl#2}\par\medbreak}

\beginsection 1.  AXLES

We assume familiarity with [\reffc]. The purpose of this manuscript
is to give a description of the program that we used to establish
[\reffc, Theorem (7.1)].

We begin by showing that it suffices to prove the equivalent of
[\reffc, Theorem (7.1)]
for parts $(K,a,b)$,  where
$a,b:V(K)\to \{5, 6,\dots, 12\}$.  While this is not really necessary, 
it makes the computer
programs slightly simpler and more elegant.  
We say that a part $(K,a,b)$ is {\it limited} if 
$a(v)\le 12$ for every $v\in V(K)$.   A {\it 
trivial limited part} is a part $(K,a,b)$ such that $b(v)=5$ and $a(v)=12$ 
for every
vertex $v$ of $K$ except the hub.  It is unique up to isomorphism.

\thm one.  Let $d=7,8,9,10,11$. If the trivial limited part with hub of 
degree
$d$ is successful, then so is the trivial part with hub of degree $d$.

\pf Let $W$ be a cartwheel with hub of degree $d$.  Let $G'$ be 
obtained from 
$G(W)$ by replacing every fan over a vertex of valency $>12$ by an 8-edge
path, and let $W'$ be the cartwheel with $G(W')=G'$ and $\gamma_{W'}(v) 
=\min\{\gamma
(v),12\}$ if $v\in V(G(W))\cap V(G')$ and $\gamma_{W'}(v)=12$ if $v\in 
V(G')-
V(G(W))$.  Then the trivial limited part of degree $d$ fits $W'$, and 
hence either a good
configuration appears in $W'$, or $N_\calp(W')\le 0$. Since every good 
configuration $K$ satisfies
$\gamma_K(v)\le 11$ for every $v\in V(K)$, and $N_\calp (W')=N_\calp(W)$ 
by condition
(iv) in the definition of a rule and the fact that $\delta (v)\in \{5, 
6, 7, 8,\infty\}$
for every rule $(G,\beta,\delta, r,s,t)$ in [\reffc, Figure 5] and
every $v\in V(G)$,  we deduce that
either a good configuration appears in $W$, or $N_\calp(W)\le 0$, as required.\qed

If $(K,a,b)$ is a part, then $K$ is, up to isomorphism, determined by 
the mappings $a,b$.  
We now make this precise.  Let $d$ be an integer.  An {\it axle of degree} 
$d$ is a pair
$A=(\ell, u)$, where $\ell, u:\{1,2,\dots, 5d\}\to\{5,6,\dots, 12\}$ such 
that
\item{(A1)} $\ell (i)\le u(i)$  for every $i=1,2\dots, 5d$, 
\item{(A2)} $\ell (i)\in\{5,6,7,8,9\}$ and $u(i)\in\{5,6,7,8,12\}$ 
for all $i=1,2,\dots, 5d$, and 

\item{(A3)} for $i=1,2,\dots, d$, if $\ell (i)\ne u(i)$, then $(\ell
(j),  u(j))=(5,12)$ for
$j=2d+i, 3d+i, 4d+i$.

\noindent We write $\ell_A=\ell$ and $u_A=u$, and put $\ell
_A(0)=u_A(0)=d$. Let $A$ be an axle of degree $d$, and let $(K,a,b)$ be
a part such that  
\item{(P1)} the hub of $K$ is 0,
\item{(P2)} the spokes of $K$ are $1,2,\dots, d$ in order,
\item{(P3)} the hats of $K$ are $d+1,d+2,\dots, 2d$ in order
(so that for $i=1,2,\dots, d-1$, $d+i$ is adjacent to $i$ and $i+1$,
and $2d$ is adjacent to $1$ and $d$), 
\item{(P4)} for $i=1,2,\dots, d$, if $k=\ell (i)=u(i)$, then 
($5\le k\le 8$ by (A2) and) the vertices of the fan over  $i$ are
$2d+i$, $3d+i,\dots, (k-4)d+i$ in order (so that $2d+i,3d+i,\dots,
(k-4)d+i,d+i$ form a path in $K$ in order; if $k=5$ there are no fan vertices),
and 
\item{(P5)} $b,a$ are the restrictions of $\ell, u$ to $V(K)$,
respectively.

\noindent In these circumstances we say that $(K,a,b)$ is the part derived 
from $A$.
It is unique up to isomorphism.

A {\it condition} is a pair $(n,m)$, where $n\in \{1,2,\dots, 5d\}$ and 
$m\in\{-8, -7, -6, -5, 6, 7, \allowbreak8, 9\}$.  We say that a
condition $(n,m)$  is 
{\it compatible} with an axle $A$ if 
\item{(C1)} $\ell_A(n)\le-m<u_A(n)$ if $m<0$,
\item{(C2)} $\ell_A(n)<m\le u_A(n)$ if $m>0$, and
\item{(C3)} either $n\le 2d$, or $n=jd+i$, where $j\in \{2,3,4\}$,
$i\in  \{1,2,\dots, d\}$
and $\ell_A(i)=u_A(i)\ge j+4$.
\smallskip

\noindent If $(n,m)$ is a condition we define $\neg(n,m)$ to be the
condition  $(n,1-m)$.
It follows immediately that $(n,m)$ is compatible with an axle if and 
only if $\neg(n,m)$ is.

Let $A$ be an axle, and let $c=(n,m)$ be a condition compatible with $A$.  
We define
$(\ell',u')$ by
$$\ell'(i)=\cases{\ell_A(i)& if $i\ne n$ or $m<0$\cr
m& otherwise}$$
$$u'(i)=\cases{u_A(i)& if $i\ne n$ or $m>0$\cr
-m& otherwise.}$$
It follows that $(\ell',u')$ is an axle; we put $A\wedge
c=(\ell',u')$.   It follows immediately that if $A$ is an axle and $c$
is a condition compatible  with $A$, then
the parts derived from $A\wedge c$ and $A\wedge (\neg c)$ are a
complementary  pair of
refinements of the part derived from $A$.  We say that an axle is {\it 
successful}
if the part derived from it is successful.  By \one\ we can restate 
[\reffc,~Theorem 
(7.1)]
as follows.  An axle $A$ of degree $d$ is {\it trivial} if $(\ell_A (i),u_A(i))=(5,12)$
for all $i=1,2,\dots, 5d$.  We denote the trivial axle of degree $d$ by 
$\uom_d$.

\thm two.  For $d=7,8,9,10,11$ the trivial axle of degree $d$ is successful.

Let $W$ be a cartwheel and $A$ an axle, both of degree $d$. We say
that $W$ is {\it compatible} with $A$ if the part derived from $A$ fits
$W$. It is easy to see that $W$ is compatible with $A$ if and only
if $W$ satisfies (P1), (P2), (P3), (P4), and $l_A(n)\le\gamma_W(n)
\le u_A(n)$ for all $n\in\{0,1,\ldots,2d\}$ and all $n$ of the
form $n=jd+i$, where $i\in\{1,2,\ldots,d\}$, $l_A(i)=u_A(i)$ and
$j=2,3,\ldots,l_A(i)-4$. We say that an axle $A$ is {\it reducible} if
for  every cartwheel
$W$ compatible with $A$ a good configuration appears in $W$.

\beginsection 2.  OUTLETS

\subsecno=2\thmno=0
Recall that a pass $P$ obeys a rule $R=(G,\beta,\delta, r,s,t)$ if $P$ 
is
isomorphic to some $(K,r,s,t)$ where $G(K)=G$ and $\beta (v)\le\gamma_K 
(v)\le
\delta (v)$ for every vertex $v\in V(G)$.  Let $h$ be the corresponding 
homeomorphism
of $\Sigma$ mapping $G(K(P))$ to $G$ and $\gamma_{K(P)}$ to $\gamma_K$.  
If
$h$ is orientation-preserving we say that $P$ {\it orientation-obeys} 
$R$; otherwise we say
that $P$ {\it anti-orientation-obeys} $R$.  If $\calr$ is a set of rules 
we write
$P\approx\calr$ to denote that $P$ orientation-obeys a member of $\calr$.

We say that a rule $R=(G,\beta,\delta, r,s,t)$ is {\it coherent} if
\item{(i)} for every cartwheel $W$ and every pass $P$ obeying $R$, if 
$P$ appears in 
$W$ in such a way that either $s(P)$ or $t(P)$ is the hub of $W$, if $v\in 
V(P)$ is
a fan of $W$ and $u\in V(G(W))$ is the unique spoke of $W$ adjacent to 
$v$ in $G(W)$,
then $u\in V(P)$ and $\beta (u')=\delta (u')$ for the corresponding vertex
$u'$ of $G$, and
\item{(ii)} if there exists a pass that both orientation-obeys $R$ and 
anti-orientation-obeys $R$,
then every pass that orientation-obeys $R$ also anti-orientation-obeys 
$R$.
\smallskip

Rules 4, 10 and 31 in [\reffc, Figure 4] are not coherent, but each can
be  split into two coherent
rules with the same net effect.  Let $\calr'$ be the set of rules obtained 
this way.  We say
that a coherent rule $R$ is {\it symmetric} if some (and hence every) 
pass that orientation-obeys
$R$ also anti-orientation-obeys $R$.  Let $\calr''$ be obtained from $\calr'$ 
by adding, 
for every non-symmetric rule $(G,\beta, \delta, r,s,t)\in\calr'$ the rule
$(G^*,\beta,\delta, r,s,t)$, where $G^*$ is isomorphic to $G$ as an abstract 
graph and
as a drawing is a ``mirror image" of $G$.  Let $\calr'''$ be obtained
from $\calr''$ by replacing the first rule by two  identical
rules of value one.
Finally, let $\calr$ be the 
set of all rules
$(G,\beta,\delta',r,s,t)$ such that $(G,\beta,\delta, r,s,t)\in\calr'''$ 
and
$\delta'(v)=\min\{\delta (v),12\}$ for every $v\in V(G)$.  Notice that
a pass may obey more than one rule in $\calr$, and that
if $(G,\beta,\delta, r,s,t)\in\calr$, then $5\le\beta\le\delta\le12$ and
$r=1$. The following holds.

\thm tone.  Let $d\ge 5$ be an integer, let $W$ be a cartwheel compatible 
with $\uom_d$,
let $w$ be the hub of $W$ and let $s$ be a spoke of $W$.  Then
$$\eqalign{&\sum (r(P):P\approx\calr, P\hbox{ appears in } W,t(P)=w, s(P)=s)\cr
&\qquad =
\sum (r(P):P\sim \calp, P \hbox{ appears in } W, t(P)=w, s(P)=s)}$$
and
$$\eqalign{&\sum (r(P):P\approx\calr, P\hbox{ appears in } W,s(P)=w, t(P)=s)\cr
&\qquad=
\sum (r(P):P\sim \calp, P \hbox{ appears in } W, s(P)=w, t(P)=s).}$$
\medskip

Now if $W$ is a cartwheel with hub $w$ and a spoke $s$ and $R\in\calr$, 
then a pass
$P\approx \{R\}$ appears in $W$ with $s(P)=s$, $t(P)=w$ if and only if 
for some vertices $v$ of $G(W),\gamma_W(v)$ is within certain bounds
determined  by $R$.
This motivates the following definition.  An {\it outlet of degree} $d$ 
is a pair $T=
(M,r)$, where $r$ is a non-zero integer, called the {\it value of} $T$, 
and $M$ is
a set $\{(p_1,\ell_1,u_1), (p_2,\ell_2,u_2),\dots, (p_n,\ell_n,u_n)\}$ 
such that
\item{(T1)} $p_1,p_2,\dots, p_n$ are integers with $1\le p_i\le 
5d$,
\item{(T2)} $\ell_1,\ell_2,\dots, \ell_n,u_1,u_2,\dots, u_n$ are
integers with $\ell_i\le u_i$, 
\item{(T3)}$\ell_i\in\{5,6,7,8,9\}$ and $u_i\in\{5,6,7,8,12\}$
for every $i=1,2,\ldots,n$, and
\item{(T4)} if $p_k=jd+i$ for some $k\in \{1,2,\dots, n\}$,
$j\in \{  2,3,4\}$
and $i\in \{1,2,\dots, d\}$, then there exists a $t\in\{1,2,\dots, n\}$ 
such that
$p_t=i$ and $l_t=u_t\ge j+4$.

\noindent We say that $T$ is {\it reduced} if $p_1,p_2,\ldots,p_n$
are pairwise distinct and $(l_i,u_i)=(5,12)$ for no $i=1,2,\ldots,n$.
We write $r(T)=r$ and $M(T)=M$.  If $i$ is an integer and
$x\in  \{0,1,\dots, d\}$
we define
$$i\oplus_d x=\cases{i+x & if $x+(i-1) \bmod d< d$\cr
i+x-d&otherwise.}$$
A {\it positioned outlet} of degree $d$ is a pair $(T,x)$, where $T$ is 
an outlet of degree $d$,
and $x\in \{1,2,\dots, d\}$.
Let $A$ be an axle of degree $d$, and let $(T,x)$ be a positioned outlet 
of degree $d$, where
$x\in \{1,2,\dots, d\}$ and 
$M(T)=\{(p_1,\ell_1, u_1), (p_2,\ell_2,u_2),\dots,
(p_n,\ell_n,u_n)\}$.
  We say that $(T,x)$ is {\it enforced} by $A$ if
$$\ell_i\le \ell_A (p_i\oplus_d (x-1))\le u_A(p_i\oplus_d (x-1))\le u_i$$
for all $i=1,2,\dots, n$.  We say that $(T,x)$ is {\it permitted by} $A$ 
if
$$u_i\ge \ell_A(p_i\oplus_d (x-1))\quad\hbox{ and } \quad
u_A(p_i\oplus_d (x-1))\ge \ell_i$$
for all $i=1,2,\dots, n$.

\thm ttwo.  For every integer $d=7,8,9,10,11$ and for every rule $R=(G,\beta,\delta,
r,s,t)\in\calr$ there exist unique reduced outlets $T$ and $T'$ such
that $r(T)=-r(T')=r$  and
for every axle $A$ and every integer $x\in \{1,2,\dots, d\}$,
\item{(i)} $(T,x)$ is enforced by $A$ if and only if for every cartwheel 
$W$ compatible
with $A$ there exists a pass $P\approx \{R\}$ appearing in $W$ with $s(P)=x$ 
and
$t(P)=0$,
\item{(ii)} $(T',x)$ is enforced by $A$ if and only if for every cartwheel 
$W$ compatible
with $A$ there exists a pass $P\approx \{R\}$ appearing in $W$ with $s(P)=0$ 
and
$t(P)=x$,
\item{(iii)} $(T,x)$ is permitted by $A$ if and only if there exist a
 cartwheel $W$ compatible
with $A$ and a pass $P\approx \{R\}$ appearing in $W$ with $s(P)=x$ and
$t(P)=0$,
\item{(iv)} $(T',x)$ is permitted by $A$ if and only if there exist a
cartwheel $W$ compatible
with $A$ and  a pass $P\approx \{R\}$ appearing 
in $W$ with $s(P)=0$ 
and
$t(P)=x$.
\smallskip

\noindent For $d=7,8,9,10,11$ let $\calt_d$ be the set of all outlets 
$T,T'$ corresponding
to rules $R\in\calr$ as in \ttwo.

Let $A$ be an axle, let $(T,x)$ be a positioned outlet, and let 
$$M(T)=\{(p_1,\ell_1,
u_1), (p_2,\ell_2,u_2),\dots, (p_n,\ell_n,u_n)\}.$$  We define
$A\wedge (T,x)$ to be the pair $A'=(\ell, u)$, where for $i=1,2,\dots, 
5d$, $l(i)$ is the least integer $l'\ge l_A(i)$ such that $l'\ge l_j$
for all $j\in\{1,2,\ldots,n\}$ with $p_j=i$, and $u(i)$ is the
largest integer $u'\le u_A(i)$ such that $u'\le u_j$ for all 
$j\in\{1,2,\ldots,n\}$ with $p_j=i$.
The following is straightforward.

\thm tthree.  Let $A$ be an axle of degree $d$, and let $(T,x)$ be a positioned 
outlet of
degree $d$.  Then $A\wedge (T,x)$ is an axle if and only if $(T,x)$ is 
permitted by $A$.

It should be noted that while $\calr$ does not depend on $d$, the corresponding 
outlets
do.  We therefore input $\calr$ in the form of a file (same for every 
$d$), and compute
the corresponding outlets of degree $d$ at the beginning of the computation.  
It is not necessary
to check correctness of this part of the program; the reader can alternatively 
verify
by inspection that the set $\calt_d$ was computed correctly.

The members of $\calr$ are stored as follows.
Let $(G,\beta,\delta, r,s,t)\in\calr$. Let us assume for convenience
that $-1\not\in V(G)$.
 We define a sequence $v_0,v_1,\dots, 
v_{16}$ such that $v_i\in V(G)\cup\{-1\}$ and every vertex of $G$
occurs in the sequence exactly once. If $u,v$
are  adjacent vertices
of $G$, let $T (u,v)$ be the vertex $w$ of $G$ such that $u,v,w$ form 
a triangle
in $G$ in clockwise order, and let $T(u,v)=-1$ if no such
vertex $w$ exists.  We define $v_0=s$, $v_1=t$, $v_2=T(v_0,v_1)$, $v_3=T(v_1,v_0)$,
$v_4=T(v_0,v_2)$, $v_5=T(v_3, v_0)$, $v_6=T(v_2,v_1)$, $v_7=T(v_1,v_3)$, 
$v_8=T(v_4,v_2)$, $v_9=T(v_3,v_5)$, $v_{10}=T(v_8,v_2)$, $v_{11}=T(v_3,v_9)$, 
$v_{12}=T(v_0,v_4)$, $v_{13}=T(v_0,v_{12})$, 
$v_{14}=T(v_5,v_0)$, $v_{15}=T(v_6,v_1)$, 
$v_{16}=T(v_{15},v_1)$. 
To input a rule we list $\beta (v_0)$, $\delta (v_0)$, $\beta (v_1)$, 
$\delta (v_1)$
and all triples $(i,\beta (v_i), \delta (v_i))$ such that 
$2\le i\le16$ and $v_i\in V(G)$.

\beginsection 3.  HUBCAPS

\subsecno=3\thmno=0
Let $d\ge 5$ be an integer.  A {\it hubcap of degree} $d$ is a collection
$((x_1,y_1,v_1),(x_2,y_2,v_2),\allowbreak
\dots,\allowbreak(x_n,y_n,v_n))$ 
of triples of integers such that 
every integer $i=1,2,\dots d$ appears in the list
$x_1,y_1,x_2,y_2,\ldots,x_n,y_n$ exactly twice.
 Let us make two remarks here. First, this definition differs somewhat 
from the one given in [\reffc]. Second, in the actual program we use
the convention that if a triple of integers appears in a hubcap twice,
it is only listed once.
Hubcaps will be used to obtain upper bounds on $N_\calp(W)$.  We now
explain  how.

Let $A$ be an axle of degree $d$, and let $x,y\in\{1,2,\dots, d\}$.  We 
put
$$L_d(A,x,y)=\max\sum_{(T,k)} r(T),$$
where the max is taken over all sets $\left\{(T_1,k_1), (T_2,k_2),\dots, 
(T_n,k_n)\right\}$
of positioned outlets with $T_i\in\calt_d$, $k_i\in \{x,y\}$ ($i=1,2,\dots, 
d$)
and such that $A'=A\wedge (T_1,k_1)\wedge (T_2,k_2)\wedge\cdots\wedge
(T_n,k_n)$ is a non-reducible axle, and the sum is over all pairs
$(T,k)$ such that $T\in\calt_d$, $k\in \{x,y\}$ and $(T,k)$ is enforced
by $A'$.  We say that  a hubcap
$H=\{(x_1,y_1,v_1),(x_2,y_2,v_2),\dots,$ $(x_n,y_n,v_n)\}$ is a {\it
hubcap} for $A$ (and that $A$ {\it has a hubcap}) if
\item{(H1)} for all $i=1,2,\dots, n$, $L_d (A,x_i,y_i)\le v_i$, and
\item{(H2)} $10(6-d)+\left\lfloor {1\over 2}\sum^n_{i=1}
v_i\right\rfloor\le  0$.

\thm thone.  Let $A$ be an axle of degree $d$ that has a hubcap.  Then 
$A$ is successful.

\pf Let $A$ be an axle of degree $d$, and let $H=\{(x_1,y_1,v_1), (x_2, 
y_2,v_2),\dots,
(x_n,y_n,v_n)\}$ be a hubcap for $A$.  Let $W$ be a cartwheel compatible 
with $A$, and
assume that no good configuration appears in $W$.  Then 
$$\eqalignno{
N_\calp(W)&=10(6-d)+\sum (r(P):
P\sim \calp, P\hbox{ appears in }W, t(P)=0)\cr
&\qquad -\sum (r(P):P\sim\calp, P \hbox{ appears
in } W, s(P)=0)\cr
&=10(6-d)+\sum (r(P):P\approx \calr, P\hbox{ appears in } W,t(P)=0)\cr
&\qquad - \sum (r(P):P\approx \calr, P \hbox{ appears in } W, s(P)=0)\cr
&= 10(6-d)+{1\over2} \sum^n_{i=1} \Bigl(\sum (r(P):P\approx \calr,
P\hbox{ appears  in }
W, t(P)=0, s(P)\in\{x_i,y_i\})\cr
&\qquad - \sum (r(P):P\approx\calr, P\hbox{ appears in } W , s(P)=0, t(P)\in
\{x_i,y_i\})\Bigr)\cr
&\le 10(6-d)+{1\over 2}\sum^n_{i=1} L_d(A,x_i,y_i)\le 10(6-d)+{1\over 
2}
\sum^n_{i=1} v_i.}$$
Since $N_\calp (W)$ is an integer, we deduce that 
$$N_\calp (W)\le 10(6-d)+\left\lfloor {1\over 2}\sum^n_{i=1} v_i\right\rfloor\le 
0,$$
as desired.\qed

We need an algorithm that given an axle $A$ and a hubcap $H$ verifies 
that $H$ is a hubcap
for $A$.  Most of that is reasonably straightforward, except for verifying 
that $L_d
(A,x,y)\le v$.  That is accomplished by a function ``CheckBound" which 
we now
describe.  Let $(T_0,z_0), (T_1,z_1),\dots,(T_{n-1},z_{n-1})$ be all the 
positioned
outlets with $T_i\in \calt_d$ and $z_i\in \{x,y\}$; the parameters of 
``CheckBound" are
integers $p\in \{0,1,\dots, n-1\}$, $s_i\in \{-1,0,1\}$ ($i=0,1,\dots, 
n-1$), $v$
and an axle $A$.  If
\item{(i)} for $i=0,1,\dots, n-1$, if $s_i=1$ then $(T_i, z_i)$ is enforced 
by $A$, and
\item{(ii)} for $i=p, p+1,\dots, n-1$ if $s_i=-1$ then $(T_i,z_i)$ is 
not permitted by $A$,

\noindent then the function ``CheckBound" verifies that either $A$ is
reducible,  or
$$\max_S \sum_j r(T_j)\le v,$$
where the max is taken over all sets $S$ such that
$$\{i\st 0\le i<n,\,s_i=1\}\subseteq S\subseteq\{i\st 0\le i<n,\,
s_i\ne-1\},$$
and  $A'=
A\wedge\bigwedge_{i\in S} (T_i,z_i)$ is an axle, and the sum is 
over all $j$ such
that $0\le j <n$ and $(T_j,z_j)$ is enforced by $A'$.  (Equivalently, the 
max
can be taken over all set $S$ such that in addition 
$r(T_i)>0$, $(T_i,z_i)$ is permitted by $A$ for every $i\in S$,
and if $(T_i,z_i)$ is enforced by $A$ then $i\in S$.)
Thus a call
to ``CheckBound" with parameters $p=0$, $s_i=0$ ($i=0,1,\dots, n-1$), 
$v$ and $A$
verifies that $L_d(A,x,y)\le v$.  The function ``CheckBound" proceeds 
in the following
steps.
\item{(1)} For every $i=p,p+1,\dots, n-1$ with $s_i=0$, if $(T_i,z_i)$ 
is enforced by
$A$ then set $s_i=1$, and if $(T_i,z_i)$ is not permitted by $A$ then 
set $s_i=-1$.
\item{(2)} Compute $f=\sum (r(T_i):0\le i<n,s_i=1)$ and 
$a=\sum (r(T_i):0\le i<n,s_i=0, r(T_i)>0)$.
\item{(3)} If $a+f\le v$ then the inequality holds.  Return.
\item{(4)} If $f>v$ test if $A$ is reducible.  If it is return, otherwise 
display an
error message and stop.
\item{(5)} For all $q=p,p+1,\dots, n-1$ with $s_q=0$ and $r(T_q)>0$ repeat
steps (6)--(10).
\item{(6)} Set $s'_i=s_i$ for $i\in \{0,1,\dots, n-1\}-\{q\}$,
$s'_q=1$, and $A'=A\wedge (T_q, z_q)$.
\item{(7)} If for some $i\in\{0,1,\dots, p-1\}$ with $s_i=-1$ the positioned 
outlet
$(T_i,z_i)$ is forced by $A'$, then go to step (9), otherwise go to step 
(8).
\item{(8)} Call ``CheckBound" recursively with arguments $q,
s'_i,v,A'$. 
\item{(9)} Set $s_q=-1$ and $a=a-r(T_q)$.
\item{(10)} If $a+f\le v$ then the inequality holds. Return.

\beginsection 4.  ASSERTIONS

\subsecno=4\thmno=0
Let $d\ge 5$ be an integer, and let $A$ be an axle of degree $d$.  We 
say that $A$ is
{\it fan-free} if $(\ell_A(i),u_A(i))=(5,12)$ for all $i=2d+1,2d+2,\dots, 
5d$.
If $A$ is fan-free we define $\tau A$ to be the axle $(\ell',u')$, where
$(\ell'(i\oplus_d 1),u'(i\oplus_d 1))=(\ell (i),u(i))$ for $i=1,2,\dots, 2d$ 
and
$(\ell'(i),u'(i))=(5,12)$ for $i=2d+1,2d+2,\dots, 5d$, and we define
$\sigma A$ to be the axle $(\ell'',u'')$, where
$$(\ell'' (i),u''(i))=\cases{
(\ell (d+1-i), u(d+1-i)) & for $i=1,2,\dots, d$\cr
(\ell (3d-i),u(3d-i)) & for $i=d+1,d+2,\dots, 2d-1$\cr
(\ell (i), u(i))& for $i=2d, 2d+1,\dots, 5d$.}$$
Thus $\tau A$ is the axle obtained from $A$ by rotating by one unit, and 
$\sigma A$ is the axle
obtained from $A$ by reflecting.  
A {\it disposition} $D$ is either $\emptyset$
 (regarded as a formal symbol informing 
us that
a certain axle is reducible), or a hubcap or a triple $(k,\ep, M)$, where 
$k$ is an
integer, $\ep\in \{0,1\}$ and $M$ is a fan-free axle.  In the first case 
we say that $D$
is a {\it reducibility disposition}, in the second case we say that $D$ 
is a {\it hubcap
disposition}, and in the third case we say that $D$ is a {\it symmetry 
disposition}.
Let $A$ be an axle, let $\calm$ be a set of axles, and let $D$ be a disposition.
We say that $D$ {\it disposes of} $A$ {\it relative to} $\calm$ if the 
following conditions hold.
\item{(i)} If $D$ is a reducibility disposition, then $A$ is reducible.
\item{(ii)} If $D$ is a hubcap disposition, then $D$ is a hubcap for $A$.
\item{(iii)} If $D=(k,\ep, M)$ is a symmetry disposition, then $M\in\calm$ 
and every
cartwheel compatible with $A$ is compatible with $\tau^k\sigma^\ep M$.
\smallskip

\noindent The following is straightforward.

\thm fone.  Let $A$ be an axle, let $\calm$ be a set of axles such that
every member of $\calm$ is successful, and let $D$ be a disposition.  
If $D$ disposes of
$A$ relative to $\calm$, then $A$ is successful.

Let $t\ge 0$ be an integer.  An {\it assertion of depth at most } $t$ 
is a sequence
$(c_1,S_1, c_2, S_2,\dots$, $c_n,S_n,D)$, where $n\ge 0$ is an integer,
$c_1,c_2,\dots, c_n$ are conditions, $S_1, S_2,\dots, S_n$ are assertions 
of depth
at most $t-1$, and $D$ is a disposition. (Thus if $t=0$, then $n=0$.) 
An {\it assertion} is an assertion of
depth at most $t$ for some integer $t\ge 0$.  
A {\it history} is a set of conditions. Let $d\ge5$ be an integer,
let $A$ be an axle of degree $d$, let $\calm$ be a 
set of fan-free axles of degree $d$,
and let $H$ be a history.  We say that an assertion $S=(c_1,S_1,c_2,S_2,\dots, 
c_n,S_n,D)$
{\it holds} for $(A,\calm, H)$ if the following two conditions are
satisfied. 
\item{(S1)} $D$ disposes of $A\wedge\bigwedge^n_{i=1} (\neg
c_i)$ relative  to $\calm\cup
\calm''$, where $\calm''$ is the set of fan-free axles of the form 
$\uom_d\wedge\bigwedge_{c\in H}
c\wedge c_i$, where $i\in \{1,2,\dots, n\}$.
\item{(S2)} For all $i=1,2,\dots, n$, $S_i$ holds for $(A',
\calm',H')$,  where
$A'=A\wedge\bigwedge_{j=1}^{i-1} (\neg c_j)\wedge c_i,$ $H'=H\cup
\{c_i\}$, and $\calm'$ consists of all members of $\calm$ and all
fan-free axles of  the form
$\uom_d\wedge\bigwedge_{c\in H}c\wedge c_j$ for $j=1,2,\dots, i-1$.
\smallskip

\noindent An assertion $S$ is a {\it presentation of degree} $d$ if $S$
holds for $(\uom_d,\emptyset,\emptyset)$.  Our proof of \two\ is based
on the following.

\thm ftwo.  If there exists a presentation of degree $d$, then $\uom_d$ 
is
successful.  

To deduce \two\ from \ftwo\ suffices to exhibit a presentation of degree 
$d$ for
every $d=7,8,9,10,11$.  Theorem \ftwo\ itself follows from the following 
more general
statement.

\thm fthree.  Let $d\ge 5$ be an integer, let $A$ be an axle of degree 
$d$, let $\calm$
be a set of fan-free axles of degree $d$, let $H$ be a history, and let $S$ be 
an assertion that holds
for $(A,\calm, H)$.  Assume that
\item{(i)} every member of $\calm$ is successful, and
\item{(ii)} for every cartwheel $W$ compatible with $\uom_d\wedge \bigwedge_{c\in 
H}c$
but not with $A$ such that $N_\calp (W)>0$, a good configuration appears in 
$W$.\ppar
\noindent Then $A$ is successful.

\pf Let $A,H,\calm, S$ be as stated.  We proceed by induction on the depth 
of $S$.
Let $S$ be of depth at most $t$, and assume that the theorem holds for 
all assertions of
depth at most $t-1$.  Let
$S=(c_1,S_1,c_2,S_2,\dots, c_n,S_n,D)$ and for $i=1,2,\dots, n$ let
$$\eqalign{
A_i&= A\wedge\bigwedge_{j<i} (\neg c_j)\wedge c_i,\cr
H_i&= H\cup \{c_i\},\cr
\calm_i&= \calm\cup \left\{\uom_d\wedge\bigwedge_{c\in
H_1}c, \uom_d\wedge\bigwedge_{c\in H_2}c,\dots,
\uom_d\wedge\bigwedge_{c\in H_{i-1}}c\right\}.}$$ We first prove the
following. \item{(1)} {\sl Let $i=1,2,\dots, n$. Then the following
statements hold.} \itemitem{(a)} {\sl For $j=1,2,\dots, i-1$, if a
cartwheel $W$ is compatible  with $\uom_d\wedge\bigwedge_{c\in H_j}c$
but not with $A_j$, then either $N_\calp (W)\le 0$ or a good
configuration appears in $W$.} \itemitem{(b)} $A_j$ {\sl is successful
for $j=1,2,\dots, i-1$.} \itemitem{(c)} {\sl Every member of $\calm_i$
is successful.}

\noindent We prove (1) by induction on $i$.  Let $i=1,2,\dots, n$, and 
assume that (a), (b) and (c)
hold for every $i'<i.$

To prove (a) we may assume that $i>1$, for otherwise (a) is vacuously 
true. It suffices
to prove the conclusion for $j=i-1$.  To this end let $W$ be compatible 
with 
$\uom_d\wedge\bigwedge_{c\in H_{i-1}}c$ but not with $A_{i-1}$.  If $W$ 
is 
not compatible with $A$ then (a) follows from (ii), and so we may assume 
that $W$ is
compatible with $A$.  Moreover, $W$ is compatible with $A\wedge c_{i-1}$, 
and so we deduce
that $W$ is compatible with one of $A_1, A_2,\dots, A_{i-2}$, and hence 
either
$N_\calp (W)\le 0$ or a good configuration appears in $W$ by the induction 
hypothesis that (b)
holds for $i-1$.  This proves (a).

To prove (b) it is enough to establish that $A_{i-1}$ is successful.  
Since $S_{i-1}$
holds for $(A_{i-1}, \calm_{i-1}, H_{i-1})$, since every member of $\calm_{i-1}$ 
is
successful by the induction hypothesis that (c) holds for every
$i'<i$, and
since every  cartwheel compatible
with $\uom_d\wedge\bigwedge_{c\in H_{i-1}}c$ but not with $A_{i-1}$
satisfies (ii) by (a) above, we deduce from the induction hypothesis that 
\fthree\ holds
for all assertions of depth at most $t-1$ that $A_{i-1}$ is successful,
as required for (b).

To prove (c) let $W$ be compatible with a member of $\calm_i$.  By (i) 
we may assume that
$W$ is compatible with $\uom_d\wedge\bigwedge_{c\in H_j} c$ for some
$j=1,2,\dots,  i-1$.
If $W$ is not compatible with $A_j$ the conclusion follows from (a); otherwise 
it
follows from (b).  This completes the proof of (1).

We are now ready to complete the proof of \fthree.  Every cartwheel compatible 
with $A$ is 
either compatible with $A_j$ for some $j=1,2,\dots, n$, or with $A'=A\wedge
\bigwedge^n_{i=1} (\neg c_i)$.  Since each $A_j$ is successful by (1b),
it  suffices
to show that $A'$ is successful.  Since $S$ holds for $(A,\calm, H)$ it follows 
that $D$
disposes of $A'$ relative to $\calm\cup\calm'$, where $\calm'$ is the 
set of fan-free
axles of the form $\uom_d\wedge\bigwedge_{c\in H}\wedge c_i$, where $i\in 
\{1,2,\dots, n\}$.
Every member of $\calm$ is successful by (i) and every member of $\calm'$ 
is successful
by (a) and (b), and hence $A'$ is successful by \fone, as required.\qed

\beginsection 5.  PRESENTATIONS

\subsecno=5\thmno=0
A presentation is described by means of a file, which in turn is described 
as a sequence
of lines.  Each line has a {\it level} associated with itself.  The lines 
are of two 
types -- {\it condition lines} describing conditions and {\it disposition 
lines} describing
dispositions.  The lines are numbered consecutively, starting from 2 (the 
first actual line of a file
describes the degree) so that line 2 has level 0.

Let $\sigma$ be a finite sequence of lines.  Let $\ell$ be the lowest 
level of a line
in $\sigma$, and let $\lam_1,\lam_2,\dots, \lam_n,\lam_{n+1}$ (in this 
order)
be all the lines of level $\ell$  in $\sigma$. If
\item{(i)} $\lam_1,\lam_{n+1}$ are the first and last lines in $\sigma$, 
respectively,
\item{(ii)} $\lam_1,\lam_2,\dots, \lam_n$ are condition lines
describing conditions $c_1,c_2,\dots, c_n$, respectively, and $\lam_{n+1}$
is a disposition line describing a disposition $D$, 
\item{(iii)} for all $i=1,2,\dots, n$, the sequence of lines of $\sigma$ 
strictly between
$\lam_i$ and $\lam_{i+1}$ describes an assertion $S_i$, and
\item{(iv)} the levels of any two consecutive lines in $\sigma$ differ 
by exactly 1,

\noindent then we say that $\sigma$ {\it describes the assertion}
$S=(c_1,S_1,c_2,S_2,\dots, c_n,S_n,D)$.  
We have created five files ``present7", ``present8",
``present9",  ``present10" and
``present11" that describe assertions $\calp_7,\calp_8,
\calp_9,\calp_{10},  \calp_{11}$,
respectively.  Our computer program verifies that they are presentations 
of appropriate degrees.

Before proceeding further we need to explain a close relationship between 
axles and outlets.
Let $T$ be an outlet of degree $d$ with $r(T)=1$ and 
$$M(T)=\{(p_1,\ell_1,u_1),(p_2,\ell_2,u_2),\dots,
(p_n,\ell_n,u_n)\}.$$
Let $A$ be the axle such that $l_i\le l_A(p_i)\le u_A(p_i)\le u_i$
for all $i=1,2,\ldots,n$,
and, subject to that, $l_A(j)\in\{5,6,7,8,9\}$ is minimum and
$u_A(j)\in\{5,6,7,8,12\}$ is maximum for all $j=1,2,\ldots,5d$.
We say that
$A$ is the {\it axle corresponding} to $T$, and that $T$ is an {\it
outlet  corresponding}
to $A$.
The {\it null condition} is the pair $(0,0)$.  If $A$ is an axle, then 
$A\wedge (0,0)$
is undefined, and hence, in particular, is not a fan-free axle.  

The program reads and processes lines of the presentation file in order.  
During execution it
maintains variables $\ell, A_i, c_i$, ($i=0,1,\dots, \ell$), $t$, $T_0, 
T_1,\dots, T_{t-1}$,
where $t,\ell$ are integers, $A_i$ are axles, $c_i$ are conditions or 
null conditions, and
$T_0, T_1,\dots, T_{t-1}$ are outlets.  At the beginning we set $\ell=t=0$,
$A_0=\uom_d$, $c_0=(0,0)$, and keep reading lines from the input file 
until all lines
are exhausted.  After reading a condition line at level $\ell$ describing 
a condition $c$ the 
program verifies that $c$ is a condition and that it is compatible with 
$A_\ell$, sets
$A_{\ell+1}=A_\ell \wedge c$ and $A_\ell = A_\ell \wedge (\neg c)$.  If
$B=\uom_d \wedge \bigwedge^\ell_{i=1}c_i$ is a fan-free axle it sets
$T_t$ to be  the outlet
corresponding to $B$ and increases $t$ by 1.  It sets $c_\ell=c$, $c_{\ell+1}=(0,0)$,
and increases $\ell$ by 1.

After reading a disposition line at level $\ell$ describing a disposition 
$D$ it verifies that $D$
disposes of $A_\ell$ relative to the set of axles corresponding to $T_0,T_1, 
\dots,
T_{t-1}$, and sets $t$ to be the largest integer $t'$ such that either 
$t'=0$ or $T_{t'-1}$
was added while executing a line of level $<\ell$.

If $\lam$ is an input line, let $\ell^{[\lam]}, t^{[\lam]}, A^{[\lam]}_i, 
c^{[\lam]}_i,
T^{[\lam]}_i$ denote the values of the variables $\ell, t, A_i, c_i, T_i$ 
immediately
prior to reading line $\lam$.  Let $\sigma^{[\lam]}$ be the sequence of 
lines consisting of $\lam$ and the
lines of level $\ge \ell^{[\lam]}$ immediately following $\lam$, and 
let
$S^{[\lam]}$ be the assertion described by $\sigma^{[\lam]}$.  Let $H^{[\lam]}$ 
consist of all $c^{[\lam]}_i$ ($i=0,1,\dots, \ell^{[\lam]}-1)$, and
let $\calm^{[\lam]}$ be the set of all axles corresponding to $T^{[\lam]}_0,
T^{[\lam]}_1,\dots, T^{[\lam]}_{t^{[\lam]}-1}$.  The following is an 
immediate consequence of
the description of the algorithm, and implies that the program
correctly verifies that  an input file describes
a presentation.

\thm fvone.  For every input line $\lam$ the program verifies that $S^{[\lam]}$
holds for $(A^{[\lam]}, \calm^{[\lam]}, H^{[\lam]})$.  

To complete the description we must explain how the program verifies disposition.
Hubcap dispositions were discussed in Section 3, reducibility dispositions 
are addressed
in the next section, and so it remains to explain how we verify symmetry 
dispositions.
A line $\lam$ describing a symmetry disposition contains four integers 
$k,\ep,\ell,m $, where $k\in \{0,1,\dots, d-1\}$, $\ep\in \{0,1\}$, $m$ 
is such 
that the line number $m$ has level $\ell$, and during its processing an 
outlet $T$ was added
to the list $T_1,\dots, T_t$ such that if $M$ denotes the axle corresponding 
to $T$, then
every cartwheel compatible with $A^{[\lam]}_{\ell^{[\lam]}}$ is compatible 
with
$\tau^k\sigma^\ep M$.  For $\ep=0$ the latter is equivalent to the fact
that $(T,k+1)$ is enforced by $A^{[\lam]}_{\ell^{[\lam]}}$, which is what 
we actually test for. We use a similar test for $\ep=1$.

\beginsection 6. TESTING APPEARANCE

\subsecno=6\thmno=0
We need to be able to test whether a given axle is reducible, and the 
purpose of this section
is to describe such test.  Let $A$ be an axle of degree $d$, and let
$B$ be the axle defined for $i=0,1,\ldots,5d$ by $(l_B(i),u_B(i))=
(u_A(i),u_A(i))$ if $1\le i\le d$ and $u_A(i)\le8$, and
$(l_B(i),u_B(i))=(l_A(i),u_A(i))$ otherwise.
Let $(K,a,b)$ be the part 
derived from $B$,
and let $L$ be the configuration with $G(L)=K$ and $\gamma_L(v)=a(v)$ 
for all
$v\in V(G(L))$.  We say that $L$ is the {\it skeleton of} $A$.  We say 
that $L$ is a {\it skeleton}
if it is a skeleton of some axle.  The {\it hub}, {\it spokes}, 
{\it hats} and {\it fans} of a  skeleton are defined in the 
obvious way.  Let $K,L$  be configurations.  We say that $K$ is a {\it 
subconfiguration}
of $L$ if $G(K)$ is a subdrawing of $G(L)$ and $\gamma_K$ is the restriction 
of
$\gamma_L$ to $V(G(K))$.  An {\it induced subconfiguration} is defined
analogously.  Thus a configuration $K$ appears in a cartwheel $W$ if
and only if $K$  is an induced
subconfiguration of $W$.  Let $L$ be a subconfiguration of a skeleton 
$K$.  We say that
$L$ is {\it well-positioned} in $K$ if for every spoke $v\in
V(G(K))-V(G(L))$, at least one of the two hats adjacent to $v$ does
not belong to $V(G(L))$.   We say that
an axle $A$ is {\it semi-reducible} if a good configuration is a well-positioned
induced subconfiguration of its skeleton.  If $K$ is such a good configuration, 
then it
follows that $K$ appears in every cartwheel $W$ compatible with $A$ such 
that
$\gamma_W(v)=\gamma_K(v)$ for every $v\in V(G(K))$.  

Later in this section we describe how we test semi-reducibility, but now, 
with that as a subroutine,
let us explain how we test whether $A$ is reducible. (Actually,
we only test for a sufficient condition for reducibility, but it
suffices for our purposes.) We
start by putting  $A$ on a stack,
and keep repeating the following steps.  
\item{(1)} If the stack is empty, then $A$ is reducible and we stop.  
Otherwise pop an axle,
say $B$, from the stack.
\item{(2)} Test if $B$ is semi-reducible.  If not then the test
failed; we display an error message and stop.
Otherwise let $L$ be a good configuration that is a well-positioned
induced subconfiguration of the skeleton of $B$.
\item{(3)} For every vertex $v$ of $G(L)$ such that $\ell_B(v)<u_B(v)$ 
do the following:
\itemitem{(a)} Let $u'_B(v)=u_B(v)-1$ and $u'_B(u)=u_B(u)$ for $u\ne v$.  
Then
$(\ell_B,u'_B)$ is an axle.
\itemitem{(b)} Put $(\ell_B,u'_B)$ on the stack.

We now explain how we test semi-reducibility, but before we do that
we should point out that verifying this part of the program is
not necessary, for there is an independent function ``CheckIso"
which (rather crudely) verifies from first principles that a mapping
produced by the semi-reducibility routine gives an isomorphism onto an
induced subconfiguration. The semi-reducibility algorithm itself
is very simple; however, its justification requires some effort.

We say that a configuration  $K$ has {\it radius
at most two} if there exists a vertex $v\in V(G(K))$ such that for every 
vertex
$u\in V(G(K))$ there is a path $P$ in $G(K)$ with ends $u,v$ and $|E(P)|\le 
2$.
The vertex $v$ is called a {\it center} of $K$.  The following is easy 
to check by inspection,
and is also verified by our computer program.

\thm sone.  Every good configuration has radius at most two.

Our semi-reducibility test is based on the following theorem. Let $L$
be a good configuration, and let $L_0$ be its free completion with ring $R$.
If $G(L)$ is $2$-connected let $J=G(L)$; otherwise there is a unique
vertex $v\in V(G(L))$ such that $G(L)\backslash v$ is disconnected.
Choose a neighbor $v'\in V(R)$ of $v$ in $L_0$, and let $J$ be the
subdrawing of $G(L_0)$ induced by $V(G(L))\cup\{v'\}$. Then $J$ is a
$2$-connected near-triangulation. In either case we say that $J$ is
an {\it enhancement} of $L$.

\thm stwo.  Let $L,K$ be configurations, let $L$ be good,
let $J$ be an enhancement of $L$,
let $J'$ be a $2$-connected near-triangulation with $V(J')=V(J)$ such that
$J'$ is a subdrawing of $J$ and there exists a $1$-$1$ mapping
$f:V(J')\to V(G(K))$ such that if  
$u,v,w\in V(J')$ form a triangle in $J'$ in the clockwise order, then
$f(u),f(v), f(w)$ form a triangle in $G(K)$ in the clockwise
order.  Assume further that $\gamma_L(v)=\gamma_K(f(v))$ for every
$v\in V(G(L))$.
Then a configuration $K_0$ isomorphic to $L$ is a
subconfiguration of $K$.  Moreover, if $L$  has radius at most two, if
$K$ is a skeleton of an axle of degree at least six, and if $K_0$ is
well-positioned in $K$, then $K_0$ is an induced subconfiguration of $K$.

\pf There exist near-triangulations $J_0=J',J_1,\dots, J_n=J$, all with 
vertex-set
$V(J)$ such that for $i=1,2,\dots, n$, $J_i$ is obtained from $J_{i-1}$ 
by adding an
edge $e_i$ with ends $u_i$ and $v_i$.  Then $e_i$ is incident with exactly 
one finite
triangle, say $T_i$, of $J_i$.  Let $w_i\not\in \{u_i,v_i\}$ be the third 
vertex
incident with $T_i$, and assume that the notation for
$u_i$ and $v_i$ is chosen so that $u_i,v_i,w_i$
form a triangle in this clockwise order.  We claim the following.
\item{(1)} {\sl For $i=1,2,\dots, n$, $f(u_i)$ is adjacent to $f(v_i)$ 
in $G(K)$,
and the vertices $f(u_i), f(v_i)$, $f(w_i)$ form a triangle in $G(K)$ 
in clockwise order.}

We prove (1) by induction on $i$.  Let $z_1=v_i,z_2,\dots,z_m=u_i$ be 
all the neighbors of $w_i$ in $J_{i-1}$
listed in the clockwise order in which they appear around $w_i$.
Since $J_{i-1}$ is a 2-connected near-triangulation we deduce that
$z_j,z_{j+1}, w_i$ form a triangle for every $j=1,2,\dots, m-1$, and
hence $f(z_j),  f(z_{j+1}),
f(w_i)$ form a triangle in $G(K)$ by the assumptions of \stwo\ and the 
induction
hypothesis.  Since $u_i,v_i$ are adjacent in $J_i$ we deduce that $w_i$ 
is not incident
with the infinite region of $J$, and hence $w_i\in V(G(L))$.  If $f(v_i)$ 
and $f(u_i)$ are
not adjacent in $G(K)$, then $\gamma_K(f(w_i))> d_{J_i}(w_i)$, 
and hence $\gamma_K(f(w_i))> d_{J_i}(w_i)=d_J(w_i)= \gamma_L(w_i)$, a 
contradiction.  Thus $f(v_i)$ and $f(u_i)$ are adjacent in $G(K)$, and 
claim (1) follows.

Next we claim
\item{(2)} {\sl If $u,v$ are adjacent in $J$, then $f(u), f(v)$ are adjacent 
in $G(K)$.}

To prove (2) let $u,v$ be adjacent in $J$, and let $P$ be a path in $J$ 
with 
vertex-set $u_0=u,u_1,\dots, u_k=v$ in order such that 
\item{(i)} $f(u_i)$ is adjacent to $f(u_{i-1})$ in $G(K)$ for every $i=1,2,\dots, 
k$,

\noindent and, subject to that,
\item{(ii)} $k$ is minimum.

\noindent Such a path exists, because every path $P$ with $E(P)\subseteq 
E(J')$ satisfies (i).
We claim that $k=1$.
To prove this we first notice that if $k=2$, then $f(u), f(v)$ are adjacent 
by (1), and so we
may assume that $k>2$.  Let $C$ be the circuit obtained from $P$ by adding 
the edge $u,v$.
Since $J$ is a near-triangulation we deduce that some pair of vertices 
$u',v'$ of $P$ other
than the two ends are adjacent in $J$.  Regardless of whether
$f(u'),f(v')$  are adjacent
in $G(K)$ or not we obtain a contradiction to the minimality of $k$.  
This proves our claim that
$k=1$, and hence completes the proof of (2).

  From (1) and (2) we deduce the first part of \stwo.  For the second
part  let $L$ have
radius at most two, and let $K$ be the skeleton of an axle of degree at 
least six.  By the first
part we may assume that $L$ is a well-positioned subconfiguration of $K$. 
If $L$ is not
an induced subconfiguration then some two vertices $u,v\in V(G(L))$
are  adjacent in
$G(K)$, but not in $G(L)$.  Since $L$ has radius at most two there exists 
a path $P$
in $G(L)$ with ends $u,v$ and $|E(P)|\le 4$.  Let us choose such a path 
$P$ with
$|E(P)|$ minimum.  Let $C$ be the circuit of $G(K)$ obtained from $P$ 
by adding
the edge $uv$.  Let $\Delta$ be the disk bounded by $C$ that is disjoint 
from the infinite
region of $G(K)$.  Since $|V(C)|\le 5$, $L$ is well-positioned in $K$ 
and $J$ is isomorphic
to a 2-connected subdrawing of $G(K)$, we deduce that $\Delta$ contains 
no vertex of $G(K)$
in its interior.  Since $G(K)$ is a near-triangulation we deduce by the 
minimality of $|E(P)|$
that $|V(P)|=3$; let $w$ be the interior vertex of $P$.  Since $u,v$
are  adjacent in $G(K)$, but not in $G(L)$, and $G(L)$ is an induced
subdrawing of $J$, 
 we deduce that $w$ is incident with the infinite region 
of $J$.
Since $J$ is 2-connected, $w$ is not incident with the infinite region
of $G(K)$, and $d_{G(K)}(w)=d_J(w)$.
 It follows that $\gamma_K(w)=d_{G(K)}(w)=d_J(w)<\gamma_L(w)$,
a contradiction which proves \stwo.\qed

Let $J$ be an enhancement of a  good configuration $L$.
For $v\in V(J)$ let $\xi (v)=\gamma_L(v)$ if $v\in V(G(L))$ and $\xi(v)=0$ 
otherwise.
A {\it query} for $L$ is a quadruple $(u,v,z,\xi (z))$, where $u,v,z$ 
are vertices of $J$
forming a triangle in the clockwise order.  A {\it question} for 
$L$ is a sequence
$Q=(Q_0,Q_1,\dots, Q_n)$ such that for $i=0,1,\dots, n$,
$Q_i=(u_i,v_i,z_i,\xi  (z_i))$ is
a query for $L$ such that $z_0,z_1,\dots, z_n$ are pairwise
distinct  and make up $V(G(L))$, $z_0$ and $z_1$ are adjacent in 
$G(L)$, and for $i=2,3,\dots, n$, $u_i,v_i\in
\{z_0,z_1,\dots,  z_{i-1}\}$.
If $Q$ is a question as above we denote by $J(Q)$ the subdrawing of $J$ 
consisting of 
all vertices of $J$, and those edges of $J$ that belong to at least
one  of the triangles
$u_i,v_i,z_i$.  Let $L,J,Q$ be as above, and let $K$ be a configuration.  
We say that $Q$
has a positive answer for $K$ if there exists a 1-1 mapping $f:V(J)\to 
V(G(K))$ such that
\item{(Q1)} $f(z_0)$ is adjacent to $f(z_1)$,
\item{(Q2)} $\gamma_K (f(z_i))= \xi (z_i)$ for all $i=0,1,\dots, n$
with  $\xi(z_i)
>0$,
\item{(Q3)} $f(u_i), f(v_i), f(z_i)$ form a triangle in $G(K)$ in the 
clockwise
order for all $i=2,3,\dots, n$.

\noindent From \stwo\ we deduce

\thm sthree.  Let $L$ be a good configuration, let $Q$ be a question 
for $L$, and
let $A$ be an axle with skeleton $K$.  If $Q$ has a positive answer
for  $K$, then a
configuration $K_0$ isomorphic to $L$ is a subconfiguration of $K$.  If 
$K_0$ is well-positioned, 
then $K_0$ is an induced subconfiguration of $K$, and hence $A$ is semi-reducible.  

Let $L$ be a good configuration, and let $Q=(u,v,z,d)$ be a query for 
$L$.  We define
$Q^*$ to be $(v,u,z,d)$.  If $Q=(Q_0,Q_1,\dots, Q_n)$ 
is a question
for $L$ we define $Q^*$, its {\it reflection}, to be
$(Q_0,Q_1,Q^*_2,  Q^*_3,
\dots, Q^*_n)$. Theorem \sthree\ does have a converse, the following.  
We omit a proof, because the result is not needed.

\thm sfour.  Let $L$ be a good configuration, let $Q$ be a question for 
$L$, let $A$
be an axle with skeleton $K$, and assume that $L$ appears in $K$.  Then 
$Q$ or $Q^*$ 
has a positive answer for $K$.

To test semi-reducibility we first compute, for every good configuration 
$L$, a question for $L$.  Then given an axle $A$ we check if $Q$ or
$Q^*$ has  a positive
answer for the skeleton $K$ of $A$.  If not, then we stop.  Otherwise 
we compute
$K_0$ as in \sthree\ and check whether it is well-positioned in $K$.

\medskip
\centerline{\bf Acknowledgment}

\noindent We would like to express our thanks to Christopher Carl Heckman
for carefully reading the manuscript and the program itself,  for
providing several useful comments, and for writing a Pascal version
of the program. Thanks also to one of the referees of [1] for pointing
out a missing assumption in an earlier version of (6.2).

\beginsection REFERENCES


\item{\reffc.}N. Robertson, D. P. Sanders, P. D. Seymour and
R. Thomas, The Four-Colour Theorem, to appear in 
{\it J.\ Combin.\ Theory Ser.\ B}.

\end

\thm stwo.  Let $L,K$ be configurations and let $J,J'$ be $2$-connected 
near-triangulations
such that $J'$ is a subdrawing of $J$, $G(L)$ is an induced
subdrawing of $J$, 
$V(J)=V(J')$, and every vertex 
in $V(J)-V(G(L))$
is incident with the infinite region of $J$.  Let $f:V(J')\to V(G(K))$ 
be a $1$-$1$ mapping
such that if $u,v,w\in V(J')$ form a triangle in $J'$ in the clockwise 
order, then
$f(u),f(v), f(w)$ form a triangle in $G(K)$ in the clockwise
order.  Assume further that $\gamma_L(v)=\gamma_K(f(v))$ for every
$v\in V(G(L))$.
Then a configuration $K_0$ isomorphic to $L$ is a
subconfiguration of $K$.  Moreover, if $L$  has radius at most two, if
$K$ is a skeleton of an axle of degree at least six, if and $K_0$ is
well-positioned in $K$, then $K_0$ is an induced subconfiguration of $K$.